# Recent developments towards optimality in multiple hypothesis testing

**Juliet Popper Shaffer**[1]

*University of California*

**Abstract:** There are many different notions of optimality even in testing a single hypothesis. In the multiple testing area, the number of possibilities is very much greater. The paper first will describe multiplicity issues that arise in tests involving a single parameter, and will describe a new optimality result in that context. Although the example given is of minimal practical importance, it illustrates the crucial dependence of optimality on the precise specification of the testing problem. The paper then will discuss the types of expanded optimality criteria that are being considered when hypotheses involve multiple parameters, will note a few new optimality results, and will give selected theoretical references relevant to optimality considerations under these expanded criteria.

## 1. Introduction

There are many notions of optimality in testing a single hypothesis, and many more in testing multiple hypotheses. In this paper, consideration will be limited to cases in which there are a finite number of individual hypotheses, each of which ascribes a specific value to a single parameter in a parametric model, except for a small but important extension: consideration of directional hypothesis-pairs concerning single parameters, as described below. Furthermore, only procedures for continuous random variables will be considered, since if randomization is ruled out, multiple tests can always be improved by taking discreteness of random variables into consideration, and these considerations are somewhat peripheral to the main issues to be addressed.

The paper will begin by considering a single hypothesis or directional hypothesis-pair, where some of the optimality issues that arise can be illustrated in a simple situation. Multiple hypotheses will be treated subsequently. Two previous reviews of optimal results in multiple testing are Hochberg and Tamhane [28] and Shaffer [58]. The former includes results in confidence interval estimation while the latter is restricted to hypothesis testing.

## 2. Tests involving a single parameter

Two conventional types of hypotheses concerning a single parameter are

$$H : \theta \leq 0 \ \ vs. \ A : \theta > 0, \tag{2.1}$$

[1]Department of Statistics, 367 Evans Hall # 3860, Berkeley, CA 94720-3860, e-mail: shaffer@stat.berkeley.edu

which will be referred to as a one-sided hypothesis, with the corresponding tests being referred to as one-sided tests, and

$$H : \theta = 0 \ \ vs, \ A : \theta \neq 0, \tag{2.2}$$

which will be referred to as a two-sided, or nondirectional hypothesis, with the corresponding tests being referred to as nondirectional tests. A variant of (2.1) is

$$H : \theta = 0 \ \ vs. \ A : \theta > 0, \tag{2.3}$$

which may be appropriate when the reverse inequality is considered virtually impossible. Optimality considerations in these tests require specification of optimality criteria and restrictions on the procedures to be considered. While often the distinction between (2.1) and (2.3) is unimportant, it leads to different results in some cases. See, for example, Cohen and Sackrowitz [14] where optimality results require (2.3) and Lehmann, Romano and Shaffer [41], where they require (2.1).

Optimality criteria involve consideration of two types of error: Given a hypothesis $H$, Type I error (rejecting $H|H$ true) and Type II error ("accepting" $H|H$ false), where the term "accepting" has various interpretations. The reverse of Type I error (accepting $H|H$ true) and Type II error (rejecting $H|H$ false, or power) are unnecessary to consider in the one-parameter case but must be considered when multiple parameters are involved and there may be both true and false hypotheses.

Experimental design often involves fixing both P(Type I error) and P(Type II error) and designing an experiment to achieve both goals. This paper will not deal with design issues; only analysis of fixed experiments will be covered.

The Neyman–Pearson approach is to minimize P(Type II error) at some specified nonnull configuration, given fixed max P(Type I error). (Alternatively, P(Type II error) can be fixed at the nonnull configuration and P(Type I error) minimized.) Lehmann [37,38] discussed the optimal choice of the Type I error rate in a Neyman-Pearson frequentist approach by specifying the losses for accepting $H$ and rejecting $H$, respectively.

In the one-sided case (2.1) and/or (2.3), it is sometimes possible to find a uniformly most powerful test, in which case no restrictions need be placed on the procedures to be considered. In the two-sided formulation (2.2), this is rarely the case. When such an ideal method cannot be found, restrictions are considered (symmetry or invariance, unbiasedness, minimaxity, maximizing local power, monotonicity, stringency, etc.) under which optimality results may be achievable. All of these possibilities remain relevant with more than one parameter, in generalized form.

A Bayes approach to (2.1) is given in Casella and Berger [12], and to (2.2) in Berger and Sellke [8]; the latter requires specification of a point mass at $\theta = 0$, and is based on the posterior probability at zero. See Berger [7] for a discussion of Bayes optimality. Other Bayes approaches are discussed in later sections.

## *2.1. Directional hypothesis-pairs*

Consider again the two-sided hypothesis (2.2). Strictly speaking, we can only either accept or reject $H$. However, in many, perhaps most, situations, if $H$ is rejected we are interested in deciding whether $\theta$ is $<$ or $>0$. In that case, there are three possible inferences or decisions: (i) $\theta > 0$, (ii) $\theta = 0$, or (iii) $\theta < 0$, where the decision (ii) is sometimes interpreted as uncertainty about $\theta$. An alternative formulation as a pair of hypotheses can be useful:

$$\begin{aligned} H_1: \theta \leq 0 \ \ vs. \ A_1: \theta > 0 \\ H_2: \theta \geq 0 \ \ vs. \ A_2: \theta < 0 \end{aligned} \tag{2.4}$$

where the sum of the rejection probabilities of the pair of tests if $\theta = 0$ is equal to $\alpha$ (or at most $\alpha$). Formulation (2.4) will be referred to as a directional hypothesis-pair.

## 2.2. Comparison of the nondirectional and directional-pair formulations

The two-sided or non-directional formulation (2.2) is appropriate, for example, in preliminary tests of model assumptions to decide whether to treat variances as equal in testing for means. It also may be appropriate in testing genes for differential expression in a microarray experiment: Often the most important goal in that case is to discover genes with differential expression, and further laboratory work will elucidate the direction of the difference. (In fact, the most appropriate hypothesis in gene expression studies might be still more restrictive: that the two distributions are identical. Any type of variation in distribution of gene expression in different tissues or different populations could be of interest.)

The directional-pair formulation (2.4) is usually more appropriate in comparing the effectiveness of two drugs, two teaching methods, etc. Or, since there might be some interest in discovering both a difference in distributions as well as the direction of the average difference or other specific distribution characteristic, some optimal method for achieving a mix of these goals might be of interest. A decision-theoretic formulation could be developed for such situations, but they do not seem to have been considered in the literature. The possible use of unequal-probability tails is relevant here (Braver [11], Mantel [44]), although these authors proposed unequal-tail use as a way of compromising between a one-sided test procedure (2.1) and a two-sided procedure (2.2).

Note that (2.4) is a multiple testing problem. It has a special feature: only one of the hypotheses can be false, and no reasonable test will reject more than one. Thus, in formulation (2.4), there are three possible types of errors:

Type I: Rejecting either $H_1$ or $H_2$ when both are true.

Type II: Accepting both $H_1$ and $H_2$ when one is false.

Type III: Rejecting $H_1$ when $H_2$ is false or rejecting $H_2$ when $H_1$ is false; i.e. rejecting $\theta = 0$, but making the wrong directional inference.

If it does not matter what conclusion is reached in (2.4) when $\theta = 0$, only Type III errors would be considered.

Shaffer [57] enumerated several different approaches to the formulation of the directional pair, variations on (2.4), and considered different criteria as they relate to these approaches. Shaffer [58] compared the three-decision and the directional hypothesis-pair formulations, noting that each was useful in suggesting analytical approaches.

Lehmann [35,37,38], Kaiser [32] and others considered the directional formulation (2.4), sometimes referring to it alternatively as a three-decision problem. Bahadur [1] treated it as deciding $\theta < 0$, $\theta > 0$, or reserving judgment. Other references are given in Finner [24].

In decision-theoretic approaches, losses can be defined as 0 for the correct decision and 1 for the incorrect decision, or as different for Type I, Type II, and Type

III errors (Lehmann [37,38]), or as proportional to deviations from zero (magnitude of Type III errors) as in Duncan's [17] Bayesian pairwise comparison method. Duncan's approach is applicable also if (2.4) is modified to eliminate the equal sign from at least one of the two elements of the pair, so that no assumption of a point mass at zero is necessary, as it is in the Berger and Sellke [8] approach to (2.1), referred to previously.

Power can be defined for the hypothesis-pair as the probability of rejecting a false hypothesis. With this definition, power excludes Type III errors. Assume a test procedure in which the probability of no errors is $1 - \alpha$. The change from a nondirectional test to a directional test-pair makes a big difference in the performance at the origin, where power changes from $\alpha$ to $\alpha/2$ in an equal-tails test under mild regularity conditions. However, in most situations it has little effect on test power where the power is reasonably large, since typically the probability of Type III errors decreases rapidly with an increase in nondirectional power.

Another simple consequence is that this reformulation leads to a rationale for using equal-tails tests in asymmetric situations. A nondirectional test is unbiased if the probability of rejecting a true hypothesis is smaller than the probability of rejecting a false hypothesis. The term "bidirectional unbiased" is used in Shaffer [54] to refer to a test procedure for (2.4) in which the probability of making the wrong directional decision is smaller than the probability of making the correct directional decision. That typically requires an equal-tails test, which might not maximize power under various criteria given the formulation (2.2).

It would seem that except for this result, which can affect only the division between the tails, usually minimally, the best test procedures for (2.2) and (2.4) should be equivalent, except that in (2.2) the absolute value of a test statistic is sometimes sufficient for acceptance or rejection, whereas in (2.4) the signed value is always needed to determine the direction. However, it is possible to contrive situations in which the optimal test procedures under the two formulations are diametrically opposed, as is demonstrated in the extension of an example from Lehmann [39], described below.

### 2.3. An example of diametrically different optimal properties under directional and nondirectional formulations

Lehmann [39] contends that tests based on average likelihood are superior to tests based on maximum likelihood, and describes qualitatively a situation in which the best symmetric test based on average likelihood is the most-powerful test and the best symmetric test based on maximum likelihood is the least-powerful symmetric test. Although the problem is not of practical importance, it is interesting theoretically in that it illustrates the possible divergence between test procedures based on (2.2) and on (2.4). A more specific illustration of Lehmann's example can be formulated as follows:

Suppose, for $0 < X < 1$, and $\gamma > 0$, known:

$$f_0(x) \equiv 1, \text{ i.e. } f_0 \text{ is Uniform } (0,1)$$
$$f_1(x) = (1+\gamma)x^\gamma,$$
$$f_2(x) = (1+\gamma)(1-x)^\gamma.$$

Assume a single observation and test $H : f_0(x)$ vs. $A : f_1(x)$ or $f_2(x)$.

One of the elements of the alternative is an increasing curve and the other is decreasing. Note that this is a nondirectional hypothesis, analogous to (2.2) above.

It seems reasonable to use a symmetric test, since the problem is symmetric in the two directions. If $\gamma > 1$ (convex $f_1$ and $f_2$), the maximum likelihood ratio test (MLR) and the average likelihood ratio test (ALR) coincide, and the test is the most powerful symmetric test:

Reject $H$ if $0 < x < \alpha/2$ or $1 - \alpha/2 < x < 1$, i.e. an extreme-tails test. However, if $\gamma < 1$ (concave $f_1$ and $f_2$), the most-powerful symmetric test is the ALR, different from the MLR; the ALR test is:

Reject $H$ if $.5 - \alpha/2 \leq x \leq .5 + \alpha/2$, i.e. a central test.

In this case, among tests based on symmetric $\alpha/2$ intervals, the MLR, using the two extreme $\alpha/2$ tails of the interval (0,1), is the *least* powerful symmetric test. In other words, regardless of the value of $\gamma$, the ALR is optimal, but coincides with the MLR only when $\gamma > 1$. Figure 1 gives the power curves for the central and extreme-tails tests over the range $0 \leq \gamma \leq 2$.

But suppose it is important not only to reject $H$ but also to decide in that case whether the nonnull function is the increasing one $f_1(x) = (1+\gamma)x^\gamma$, or the decreasing one $f_2(x) = (1+\gamma)(1-x)^\gamma$. Then the appropriate formulation is analogous to (2.4) above:

$$H_1 : f_0 \text{ or } f_1 \text{ vs. } A_1 : f_2$$
$$H_2 : f_0 \text{ or } f_2 \text{ vs. } A_2 : f_1.$$

If $\gamma > 1$ (convex), the most-powerful symmetric test of (2.2) (MLR and ALR) is also the most powerful symmetric test of (2.4). But if $\gamma < 1$ (concave), the most-powerful symmetric test of (2.2) (ALR) is the least-powerful while the MLR is the most-powerful symmetric test of (2.4). In both cases, the directional ALR and MLR are identical, since the alternative hypothesis consists of only a single distribution. In general, if the alternative consists of a single distribution, regardless of the dimension of the null hypothesis, the definitions of ALR and MLR coincide.

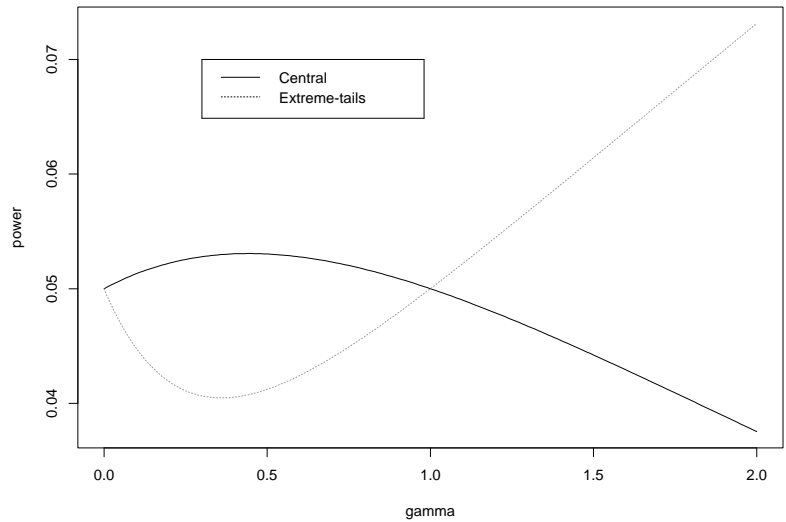

Fig 1. *Power of nondirectional central and extreme-tails tests*

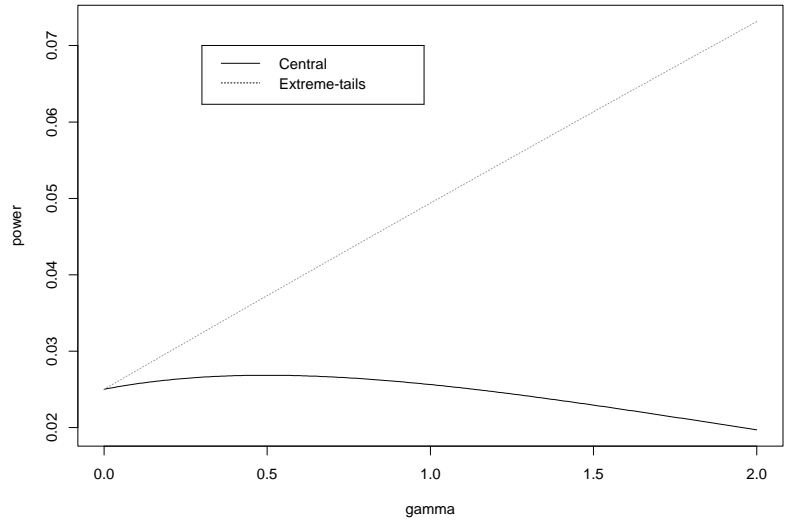

Fig 2. *Power of directional central and extreme-tails test-pairs*

Note that if $\gamma$ is unknown, but is known to be $< 1$, the terms 'most powerful' and 'least powerful' in the example can be replaced by 'uniformly most powerful' and 'uniformly least powerful', respectively.

Figure 2 gives some power curves for the directional central and extreme-value test-pairs.

Another way to look at this directional formulaton is to note that a large part of the power of the ALR under (2.2) becomes Type III error under (2.4) when $\gamma < 1$. The Type III error not only stays high for the directional test-pair based on the central proportion $\alpha$, but it actually is above the null value of $\alpha/2$ when $\gamma$ is close to zero.

Of course the situation described in this example is unrealistic in many ways, and in the usual practical situations, the best test for (2.2) and for (2.4) are identical, except for the minor tail-probability difference noted. It remains to be seen whether there are realistic situations in which the two approaches diverge as radically as in this example. Note that the difference between nondirectional and directional optimality in the example generalizes to multiple parameter situations.

## 3. Tests involving multiple parameters

In the multiparameter case, the true and false hypotheses, and the acceptance and rejection decisions, can be represented in a two by two table (Table 1). With more than one parameter, the potential number of criteria and number of restrictions on types of tests is considerably greater than in the one-parameter case. In addition, different definitions of power, and other desirable features, can be considered. This paper will describe some of these expanded possibilities. So far optimality results have been obtained for relatively few of these conditions. The set of hypotheses to be considered jointly in defining the criteria is referred to as the family. Sometimes the family includes all hypotheses to be tested in a given study, as is usually the case, for example, in a single-factor experiment comparing a limited number of treatments. Hypotheses tested in large surveys and multifactor experiments are usually divided into subsets (families) for error control. Discussion of choices for families can be found in Hochberg and Tamhane [28] and Westfall and Young [66].

TABLE 1
*True states and decisions for multiple tests*

| Number of | Number not rejected | Number rejected | |
|---|---|---|---|
| True null hypotheses | U | V | $m_0$ |
| False null hypotheses | T | S | $m_1$ |
| | m-R | R | $m$ |

All of the error, power, and other properties raise more complex issues when applied to tests of (2.2) than to tests of (2.1) or (2.3), and even more so to tests of (2.4) and its variants. With more than one parameter, in addition to these expanded possibilities, there are also more possible types of test procedures. For example, one may consider only stepwise tests, or, even more specifically, under appropriate distributional assumptions, only stepwise tests using either $t$ tests or $F$ tests or some combination. Some type of optimality might then be derived within each type. Another possibility is to derive optimal results for the sequence of probabilities to be used in a stepwise procedure without specifying the particular type of tests to be used at each stage. Optimality results may also depend on whether or not there are logical relationships among the hypotheses (for example when testing equality of all pairwise differences among a set of parameters, transitivity relationships exist). Other results are obtained under restrictions on the joint distributions of the test statistics, either independence or some restricted type of dependence. Some results are obtained under the restriction that the alternative parameters are identical.

### 3.1. Criteria for Type I error control

Control of Type I error with a one-sided or nondirectional hypothesis, or Type I and Type III error with a directional hypothesis-pair, can be generalized in many ways. Type II error, except for one definition below (viii), has usually been treated instead in terms of its obverse, power. Optimality results are available for only a small number of these error criteria, mainly under restricted conditions.

Until recent years, the generalized Type I error rates to be controlled were limited to the following three:

(i) The expected proportion of errors (true hypotheses rejected) among all hypotheses, or the maximum per-comparison error rate (PCER), defined as E(V/m). This criterion can be met by testing each hypothesis at the specified level, independent of the number of hypotheses; it essentially ignores the multiplicity issue, and will not be considered further.
(ii) The expected number of errors (true hypotheses rejected), or the maximum per-family error rate (PFER), where the family refers to the set of hypotheses being treated jointly, defined as E(V).
(iii) The maximum probability of one or more rejections of true hypotheses, or the familywise error rate (FWER), defined as Prob($V > 0$).

The criterion (iii) has been the most frequently adopted, as (i) is usually considered too liberal and (ii) too conservative when the same fixed conventional level is adopted. Within the last ten years, some additional rates have been proposed to meet new research challenges, due to the emergence of new methodologies and technologies that have resulted in tests of massive numbers of hypotheses and a concomitant desire for less strict criteria.

Although there have been situations for some time in which large numbers of hypotheses are tested, such as in large surveys, and multifactor experimental designs, these hypotheses have usually been of different types and often of an indefinite number, so that error control has been restricted to subsets of hypotheses, or families, as noted above, each usually of some limited size. Within the last 20 years, there has been an explosion of interest in testing massive numbers of well-defined hypotheses in which there is no obvious basis for division into families, such as in microarray genomic analysis, where individual hypotheses may refer to parameters of thousands of genes, to tests of coefficients in wavelet analysis, and to some types of tests in astronomy. In these cases the criterion (iii) seems to many researchers too draconian. Consequently, some new approaches to error control and power have been proposed. Although few optimal criteria have been obtained under these additional approaches, these new error criteria will be described here to indicate potential areas for optimality research.

Recently, the following error-control criteria in addition to (i)-(iii) above have been considered:

(iv) The expected proportion of falsely-rejected hypotheses among the rejected hypotheses–the false discovery rate (FDR). The proportion itself, FDP = V/R, is defined to be 0 when no hypotheses are rejected (Benjamini and Hochberg [3]; for earlier discussions of this concept see Eklund and Seeger [22], Seeger [53], Sorić [59]), so the FDR can be defined as $E(FDP|R > 0)P(R > 0)$. There are numerous publications on properties of the FDR, with more appearing continuously.

(v) The expected proportion of falsely-rejected hypotheses among the rejected hypotheses given that some are rejected ($p$-FDR) (Storey [62]), defined as $E(V/R)|R > 0$.

(vi) The maximum probability of at most $k$ errors ($k$-FWER or $g$-FWER–$g$ for generalized), given that at least $k$ hypotheses are true, $k = 0, \ldots, m$, $P(V > k)$, (Dudoit, van der Laan, and Pollard [16], Korn, Troendle, McShane, and Simon [34], Lehmann and Romano [40], van der Laan, Dudoit, and Pollard [65], Pollard and van der Laan [48]). Some results on this measure were obtained earlier by Hommel and Hoffman [31].

(vii) The maximum proportion of falsely-rejected hypotheses among those rejected (with 0/0 defined as 0), $FDP > \gamma$ (Romano and Shaikh [51] and references listed under (vi)).

(viii) The false non-discovery rate (Genovese and Wasserman [25]), the expected proportion of nonrejected but false hypotheses among the nonrejected ones, (with 0/0 defined as 0): $FNR = E[T/(m - R) \, P(m - R > 0)]$.

(ix) The vector loss functions defined by Cohen and Sackrowitz [15], discussed below and in Section 4.

Note that the above generalizations (iv), (vi), and (vii) reduce to the same value (Type I error) when a single parameter is involved, (v) equals unity so would not be appropriate for a single test, (viii) reduces to the Type II error probability, and the FRR in (ix), defined in Section 4, is equal to (ii).

The loss function approach has been generalized as either ($L_i$) the sum of the loss functions for each hypothesis, ($L_{ii}$) a 0-1 loss function in which the loss is zero only if all hypotheses are correctly classified, or ($L_{iii}$) a sum of loss functions for the FDR (*iv*) and the FNR (*viii*). (In connection with $L_{iii}$, see the discussion of Genovese and Wasserman [25] in Section 4, as well as contributions by Cheng et al [13], who also consider adjusting criteria to consider known biological results in genomic applications.) Sometimes a vector of loss functions is considered rather than a composite function when developing optimal procedures; a number of different vector approaches have been used (see the discussion of Cohen and Sackrowitz [14,15] in the section on optimality results). Many Bayes and empirical Bayes approaches involve knowing or estimating the proportion of true hypotheses, and will be discussed in that context below.

The relationships among (iv), (v), (vi) and (vii) depend partly on the variance of the number of falsely-rejected hypotheses. Owen [47] discusses previous work related to this issue, and provides a formula that takes the correlations of test statistics into account.

A contentious issue relating to generalizations (iv) to (vii) is whether rejection of hypotheses with very large *p*-values should be permitted in achieving control using these criteria, or whether some additional restrictions on individual *p*-values should be applied. For example, under (vi), $k$ hypotheses could be rejected even if the overall test was negative, regardless of the associated *p*-values. Under (iv), (v), and (vii), given a sufficient number of hypotheses rejected with FWER control at $\alpha$, additional hypotheses with arbitrarily large p-values can be rejected. Tukey, in a personal oral communication, suggested, in connection with (iv), that hypotheses with individual *p*-values greater than $\alpha$ might be excluded from rejection. This might be too restrictive, especially under (vi) and (vii), but some restrictions might be desirable. For example, if $\alpha = .05$, it has been suggested that hypotheses with $p \leq \alpha*$ might be considered for rejection, with $\alpha*$ possibly as large as 0.5 (Benjamini and Hochberg [4]). In some cases, it might be desirable to require nonincreasing individual rejection probabilities as $m$ increases (with $m \geq k$ in (vi) and (vii)), which would imply Tukey's suggestion. Even the original Benjamini and Hochberg [3] FDR-controlling procedure violates this latter condition, as shown in an example in Holland and Cheung [29], who note that the adaptive method in Benjamini and Hochberg [4] violates even the Tukey suggestion. Consideration of these restrictions on error is very recent, and this issue has not yet been addressed in any serious way in the literature.

### *3.2. Generalizations of power, and other desirable properties*

The most common generalizations of power are:

(a) probability of at least one rejection of a false hypothesis, (b) probability of rejecting all false hypotheses, (c) probability of rejecting a particular false hypothesis, and (d) average probability of rejecting false hypotheses. (The first three were initially defined in the paired-comparison situation by Ramsey [49], who called them any-pair power, all-pairs power, and per-pair power, respectively.) Generalizations (b) and (c) can be further extended to (e) the probability of rejecting more than $k$ false hypotheses, $k = 0, \ldots, m$. Generalization (d) is also the expected proportion of false hypotheses rejected.

Two other desirable properties that have received limited attention in the literature are:

(f) complexity of the decisions. Shaffer [55] suggested the desirability, when comparing parameters, of having procedures that are close to partitioning the parameters into groups. Results that are partitions would have zero complexity; Shaffer suggested a quantitative criterion for the distance from this ideal.

(g) familywise robustness (Holland and Cheung [29]). Because the decisions on definitions of families are subjective and often difficult to make, Holland and Cheung suggested the desirability of procedures that are less sensitive to family choice, and developed some measures of this criterion.

## 4. Optimality results

Some optimality results under error protection criteria (i) to (iv) and under Bayesian decision-theoretic approaches were reviewed in Hochberg and Tamhane [28] and Shaffer [58]. The earlier results and some recent extensions will be reviewed below, and a few results under (vi) and (viii) will be noted. Criteria (iv) and (v) are asymptotically equivalent when there are some false hypotheses, under mild assumptions (Storey, Taylor and Siegmund [63]).

Under (ii), optimality results with additive loss functions were obtained by Lehmann [37,38], Spjøtvoll [60], and Bohrer [10], and are described in Shaffer [58] and Hochberg and Tahmane [28]. Lehmann [37,38] derived optimality under hypothesis-formulation (2.1) for each hypothesis, Spjøtvoll [60] under hypothesis-formulation (2.2), and Bohrer [10] under hypothesis-formulation (2.4), modified to remove the equality sign from one member of the pair.

Duncan [17] developed a Bayesian decision-theoretic procedure with additive loss functions under the hypothesis-formulation (2.4), applied to testing all pairwise differences between means based on normally-distributed observations, and assuming the true means are normally-distributed as well, so that the probability of two means being equal is zero. In contrast to Lehmann [37,38] and Spjøtvoll [60], Duncan uses loss functions that depend on the magnitudes of the true differences when the pair (2.4) are accepted or the wrong member of the pair (2.4) is rejected. Duncan [17] also considered an empirical Bayes version of his procedure in which the variance of the distribution of the true means is estimated from the data, pointing out that the results were almost the same as in the known-variance case when $m \geq 15$. For detailed descriptions of these decision-theoretic procedures of Lehmann, Spjøtvoll, Bohrer, and Duncan, see Hochberg and Tamhane [28].

Combining (iv) and (viii) as in $L_{iii}$, Genovese and Wasserman [25] consider an additive risk function combining FDR and FNR and obtain some optimality results, both finite-sample and asymptotic. If the risk $\delta_i$ for $H_i$ is defined as 0 when the correct decision is made (either acceptance or rejection) and 1 otherwise, they define the classification risk as

$$R_m = \frac{1}{m} E(\sum_{i=1}^{m} |\delta_i - \hat{\delta}_i|), \qquad (4.1)$$

equivalent to the average fraction of errors in both directions. They derive asymptotic values for $R_m$ given various procedures and compare them under different conditions. They also consider the loss functions $FNR + \lambda FDR$ for arbitrary $\lambda$ and derive both finite-sample and asymptotic expressions for minimum-risk procedures based on $p$-values. Further results combining (iv) and (viii) are obtained in Cheng et al [13].

Cohen and Sackrowitz [14,15] consider the one-sided formulations (2.1) and (2.3) and treat the multiple situation as a $2^m$ finite action problem. They assume a multivariate normal distribution for the $m$ test statistics with a known covariance matrix of the intraclass type (equal variances, equal covariances), so the test statistics are exchangeable. They consider both additive loss functions (1 for a Type I error and an arbitrary value $b$ for a Type II error, added over the set of hypothesis tests), the m-vector of loss functions for the m tests, and a 2-vector treating Type I and Type II errors separately, labeling the two components false rejection rate (FRR) and false acceptance rate (FAR). They investigate single-step and stepwise procedures from the points of view of admissibility, Bayes, and limits of Bayes procedures. Among a series of Bayes and decision-theoretic results, they show admissibility of single-stage and stepdown procedures, with inadmissibility of stepup procedures, in contrast to the results of Lehmann, Romano and Shaffer [41], described below, which demonstrate optimality of stepup procedures under a different loss structure.

Under error criterion (iii), early results on means are described in Shaffer [58] with references to relevant literature. Lehmann and Shaffer [42], considering multiple range tests for comparing means, found the optimal set of critical values assuming it was desirable to maximize the minimum probabilities for distinguishing among pairs of means, which implies maximizing the probabilities for comparing adjacent means. Finner [23] noted that this optimality criterion was not totally compelling, and found the optimal set under the assumption that one would want to maximize the probability of rejecting the largest range, then the next-largest, etc. He compared the resulting maximax method to the Lehmann and Shaffer [42] maximin method.

Shaffer [56] modified the empirical Bayes version of Duncan [17], decribed above, to provide control of (iii). Recently, Lewis and Thayer [43] adopted the Bayesian assumption of a normal distribution of true means as in Duncan [17]and Shaffer [56] for testing the equality of all pairwise differences among means. However, they modified the loss functions to a loss of 1 for an incorrect directional decision and $\alpha$ for accepting both members of the hypothesis-pair (2.4). False discoveries are rejections of the wrong member of the pair (2.4) while true discoveries are rejections of the correct member of the pair. Thus, Lewis and Thayer control what they call the directional FDR (DFDR). Under their Bayesian and loss-function assumptions, and adding the loss functions over the tests, they prove that the DFDR of the minimum Bayes-risk rule is $\leq \alpha$. They also consider an empirical Bayes variation of their method, and point out that their results provide theoretical support for an empirical finding of Shaffer [56], which demonstrated similar error properties for her modification of the Duncan [17] approach to provide control of the FWER and the Benjamini and Hochberg [3] FDR-controlling procedure. Lewis and Thayer [43] point out that the empirical Bayes version of their assumptions can be alternatively regarded as a random-effects frequentist formulation.

Recent stepwise methods have been based on Holm's [30] sequentially- rejective procedure for control of (iii). Initial results relevant to that approach were obtained by Lehmann [36], in testing one-sided hypotheses. Some recent results in Lehmann, Romano, and Shaffer [41] show optimality of stepwise procedures in testing one-sided hypotheses for controlling (iii) when the criterion is maximizing the minimum power in various ways, generalizing Lehmann's [36] results. Briefly, if rejection of at least $i$ hypotheses are ordered in importance from $i = 1$ (most important) to $i = m$, a generalized Holm stepdown procedure is shown to be optimal, while if these rejections are ordered in importance from $i = m$ (most important) to $i = 1$, a stepup procedure generalizing Hochberg [27] is shown to be optimal.

Most of the recent literature in multiple comparisons relates to improving existing methods, rather than obtaining optimal methods, although of course such improvements indicate the directions in which optimality might be achieved. The next section is a necessarily brief and selective overview of some of this literature.

## 5. Literature on improvement of multiple comparison procedures

### 5.1. Estimating the number of true null hypotheses

Under most types of error control, if the number $m_0$ of true hypotheses $H$ were known, improved procedures could be based on this knowledge. In fact, sometimes such knowledge could be important in its own right, for example in microarray analysis, where it might be of interest to estimate the number of genes differentially expressed under different conditions, or in an astronomy problem (Meinshausen and Rice [45]) in which that is the only quantity of interest.

Under error criteria (ii) and (iii), the Bonferroni method could be improved in power by carrying it out at level $\alpha/m_0$ instead of $\alpha/m$. FDR control with independent test statistics using the Benjamini and Hochberg [3] method is exactly equal to $\pi_0 \alpha$, where $\pi_0 = m_0/m$ is the proportion of true hypotheses, (Benjamini and Hochberg [3]), so their FDR-controlling method described in that paper could be made more powerful by multiplying the criterion $p$-values at each stage by $m/m_0$. The method has been proved to be conservative under some but not all types of dependence. A modified method making use of $m_0$ guarantees control of the FDR at the specified level under all types of test statistic dependence (Benjamini and Yekutieli [6]).

Much recent effort has been directed towards obtaining good estimates of $\pi_0$, either for an interest in this quantity itself, or because then these improved methods and other more recent methods, including some single-step methods, could be used at level $\alpha* = \alpha/\pi_0$. There are many recent papers comparing estimates of $\pi_0$, but few optimality results are available at present. Some recent relatively theoretical references are Black [9], Genovese and Wasserman [26], Storey, Taylor, and Siegmund [63], Meinshausen and Rice [45], and Reiner, Yekutieli and Benjamini [50].

Storey, Taylor, and Siegmund [63] use empirical process theory to investigate proposed procedures for FDR control. The original Benjamini and Hochberg [3] procedure is a stepup procedure, using the ordered $p$-values, while the procedures proposed in Storey [61] and others are single-step procedures in which all $p$-values less than a criterion $t$ are rejected. Based on the notation in Table 1, Storey, Taylor and Siegmund [63] define the empirical processes

$$\begin{aligned} V(t) &= \#(null\ p_i : p_i \leq t) \\ S(t) &= \#(alternative\ p_i : p_i \leq t) \\ R(t) &= V(t) + S(t) = \#(p_i : p_i \leq t). \end{aligned} \quad (5.1)$$

They use empirical process theory to prove both finite-sample and asymptotic control of FDR for the Benjamini and Hochberg [3] procedure and the most conservative Storey [61] procedure, and also for new proposed procedures that involve estimation of $\pi_0$ under both independence and some forms of positive dependence.

Benjamini, Krieger, and Yekutieli [5] develop two-stage and multistage adaptive methods, and study the two-stage method analytically. That method provides an estimate of $\pi_0$ at the first stage and takes the uncertainty about the estimate into account in modifying the second stage. It is proved to guarantee FDR control at the specified level. Based on extensive simulation results the methods proposed in Storey, Taylor and Siegmund [63] perform best when test statistics are independent, while the Benjamini, Krieger and Yekutieli [5] two-stage adaptive method appears to be the only proposed method (to this time) based on estimating $m_0$ that controls the FDR under the conditions of high positive dependence that are sufficient for FDR control using the original Benjamini and Hochberg [3] FDR procedure.

### 5.2. Resampling methods

In general, under any criterion, if appropriate aspects of joint distributions of test statistics were known (e.g. their covariance matrices), procedures based on those distributions could achieve greater power with the same error control than procedures ensuring error control but not based on such knowledge. Resampling methods are being intensively investigated from this point of view. Permutation methods, when applicable, can provide exact error control under criterion (iii) (Westfall and Young [66]) and some bootstrap methods have been shown to provide asymptotic error control, with the possibility of finding asymptotically optimal methods under such control (Dudoit, van der Laan and Pollard [16], Korn, Troendle, McShane and Simon [34], Lehmann and Romano [40], van der Laan, Dudoit and Pollard [65], Pollard and van der Laan [48], Romano and Wolf [52]).

Since the asymptotic methods are based on the assumption of large sample sizes relative to the number of tests, it is an open question how well they apply in cases of massive numbers of hypotheses in which the sample size is considerably smaller than $m$, and therefore how relevant any asymptotic optimal properties would be in these contexts. Some recent references in this area are Troendle, Korn, and McShane [64], Bang and Young [2]), and Muller et al [46], the latter in the context of a Bayesian decision-theoretic model.

### 5.3. Empirical Bayes procedures

The Bayes procedure of Berger and Sellke [8] referred to in the section on a single parameter, testing (2.1) or (2.2), requires an assumption of the prior probability that the hypothesis is true. With large numbers of hypotheses, the procedure can be replaced by empirical Bayes procedures based on estimates of this prior probability by estimating the proportion of true hypotheses. These, as well as estimates of other aspects of the prior distributions of the test statistics corresponding to true and false hypotheses, are obtained in many cases by resampling methods. Some of the references in the two preceding subsections are relevant here; see also Efron [18], and Efron and Tibshirani [21]; the latter compares an empirical Bayes method with the FDR-controlling method of Benjamini and Hochberg [3]. Kendziorski et al [33] use an empirical Bayes hierarchical mixture model with stronger parametric assumptions, enabling them to estimate the relevant parameters by log likelihood methods rather than resampling.

For an unusual approach to the choice of null hypothesis, see Efron [19], who suggests that an alternative null hypothesis distribution, based on an empirically-determined "central" value, should be used in some situations to determine "interesting" – as opposed to "significant" – results. For a novel combination of empirical Bayes hypothesis testing and estimation, related to Duncan's [17] emphasis on the magnitude of null hypothesis departure, see Efron [20].

## 6. Summary

Both one-sided and two-sided tests referring to a single parameter are considered. A two-sided test referring to a single parameter becomes multiple inference when the hypothesis that the parameter $\theta$ is equal to a fixed value $\theta_0$ is reformulated as the directional hypothesis-pair (i) $\theta \leq \theta_0$ and (ii) $\theta \geq \theta_0$, a more appropriate formulation when directional inference is desired. In the first part of the paper, it is shown that optimality results in the case of a single nondirectional hypothesis can be diametrically opposite to directional optimality results. In fact, a procedure that is uniformly most powerful under the nondirectional formulation can be uniformly least powerful under the directional hypothesis-pair formulation, and vice versa.

The second part of the paper sketches the many different formulations of error rates, power, and classes of procedures when there are multiple parameters. Some of these have been utilized for many years, and some are relatively new, stimulated by the increasing number of areas in which massive sets of hypotheses are being tested. There are relatively few optimality results in multiple comparisons in general, and still fewer when these newer criteria are utilized, so there is great potential for optimality research in this area. Many existing optimality results are described in Hochberg and Tamhane [28] and Shaffer [58]). These are sketched briefly here, and some further relevant references are provided.